\documentclass[twoside]{article}
\input{spsp10.sty}

\begin{document}
\renewcommand{\baselinestretch}{1}

\thispagestyle{empty}
\fancyhead[LE]{Reproducing kernel method for PDE constrained optimization}

\fancyhead[RO]{Darehmiraki}

\begin{center}
\vspace*{0.1cm}
{\Large \bf Reproducing kernel method for PDE constrained optimization}\\
{
Majid Darehmiraki {\\ darehmiraki@gmail.com }\\
Department of Mathematics, Khatam Al-Anba university of technology, Behbahan, Iran}

\end{center}
\noindent{\bf Abstract:}
This paper presents  reproducing kernel Hilbert spaces method (RKHSM) to obtain the numerical
solutions for PDE constrained optimization problem. The analytical solution is shown in a series form in the reproducing
kernel space and the approximate solution is constructed by truncating the series. Convergence
analysis of presented method is discussed.  Several test problems are employed and results of numerical experiments are presented. The obtained results confirm the acceptable accuracy of the proposed methods.    

\section{Introduction}
The parabolic optimal control problems can be written as follows:
We consider a bounded Lipschitz domain $\Omega\subset\mathbb{R}^N,N\geq1$ with boundary $ \Gamma $. The domain $ \Omega $  stands for a spatial domain that is to be heated in the fixed time interval $ [0,T] $. The heating is done by a controlled heat source of density $ u:\Omega\times [0,T]\longrightarrow\mathbb{R} $. \\ 
By $ y(x,t) $ we mean the temperature in the point $ x\in\Omega $ at time $ t\in [0,T] $ while $ y_0(x) $ is the known temperature at the initial time $ t=0 $. We assume that, for given $ u $, the temperature $ y $ is obtained by the solution of the following linear heat equation
\begin{align}
\begin{cases}
-\dfrac{\partial y}{\partial t}(x,t)+\triangle y(x,t)\, & =\,u(x,t)\,\,\,  in\,\,  Q:=\Omega\times (0,T)\\
\qquad y(x,t)\, & =\,h(x,t)\,\,\, on\,\,  \Sigma :=\Gamma\times (0,T)\\
\qquad  y(x,0)\, & =\,y_0(x)\,\,\,  in\,\,  \Omega
\end{cases}
\end{align}
The Dirichlet boundary condition says that the temperature $ y $ at the boundary $ \Gamma $ is zero at any time. The function $ u $ is the control function, while $ y $ is called the associated state; the partial differential equation (1.1) is said to be the state equation, and is a linear parabolic equation. In this formulation, we have tacitly assumed that to each control $ u $ there exist a unique state $ y $.\\
There are many possible forms of the cost functional to be minimized. Perhaps the most common model, and the one we examine in this article, is of the form 
$$(P)\,\,\,\,\min\limits_{u\in U_{ad}}\,\,J(u):=\dfrac{1}{2	}\int_0^T\int_{\Omega} (y(x,t)-y^d(x,t))^2d\Omega dt+\dfrac{\nu}{2	}\int_0^T\int_{\Omega}u^2(x,t)d\Omega dt.$$ 
The quantity $ \nu >0 $ is the Tikhonov regularization parameter, which determines to what extent one wishes to achieve realization of the desired state and minimization of the control. Now we have an optimal control problem consists of the above cost functional and the linear parabolic equation (1.1). The integral funcional above is convex with respect to $ y $ . The problem $ P $ is in general convex, because the equations (1.1) is linear \cite{Borzi}. 

\subsection{The literature review}
Distributed optimal control has multitude applications in science and engineering. For instance, distributed optimal control problems arise in such diverse areas as aerodynamic, mathematical finance, medicine, and environmental engineering. As the computational capacity increases and optimization techniques become more advanced offer the possibility to solve optimization problems easier and faster. \\The distributed optimal control problems are generally difficult to solve and their exact solutions are difficult to obtain, therefore, some various approximate methods have recently been developed such as radial basic function method \cite{Rad}, \cite{Pearson}, fictitious domain method \cite{Eppler}, proper orthogonal decomposition \cite{Ravindran}, interior point method \cite{Weiser}, Newton method \cite{Laumen}, domain decomposition method \cite{Benamou}, and the Variational Iteration Method \cite{Akkouche}. \\The numerical methods used to find the optimal control of parabolic distributed parameter systems
have been presented by \cite{Sage, Mahapatra, Wang, Horng, Chang}. Sage and White \cite{Sage} used a finite difference technique,
Mahapatra \cite{Mahapatra} derived a piecewise continuous solution using Walsh functions and Wang and Chang
\cite{Wang} transformed the optimal control problem into a two-point boundary value problem and obtained
the solution using shifted Legendre polynomials. Horng and Chou \cite{Horng} reduced the optimal control of
a distributed parameter system into the optimal control of a linear time-invariant lumped parameter
system; furthermore, they derived the integral of the cross-product of two shifted Chebyshev vectors
to find the solution. Chang and Yang \cite{Chang} transformed the optimal control problem into a two point
boundary value problem; they also derived the operational matrix for the integration of the generalized
orthogonal polynomials and obtained the optimal control using the Taylor series and several kinds of
orthogonal polynomials, by employing only the cross-product of two shifted Legendre vectors. Razzaghi
and Arabshahi \cite{Razzaghi} transformed the optimal control problem into a two point boundary value problem
and adopted an approach using the Taylor series. Rad et al. \cite{Rad} solved parabolic optimal control problem by radial basis function. For more references in this content one can see \cite{Sadek}, \cite{Kar}.\\
\subsection{The main aim of this paper}
In this study, a new iterative algorithm for solving the PDE constrained optimization problem in the reproducing kernel space is proposed. The
advantages of the approach lie in the following facts. The approximate solution of state and control functions  converges uniformly to their exact
solutions. The method is mesh free, easily implemented and capable of treating the boundary conditions. Since the
method needs no time discretization, it does not matter at what time the approximate solution is computed, from both the
elapsed CPU time and stability problem points of view. \\
The advantages of the approach lie in the following facts. The method is
mesh free, easily implemented and capable in treating various boundary conditions.
The method needs no time discretization against [7,10] and any ODE
integrator against [5, 9, 18]. Therefore there is no concern about the stability
problem and also increasing the end of time T does not increase the CPU
time.\\
This paper is arranged in the following manner, in Section 2 a brief introduction
of the reproducing kernel Hilbert spaces and several reproducing kernel spaces are represented. Section 3 present necessary optimality conditions of mentioned optimal control. The problem solving,
method implementation and verification of convergence of the approximate
solution to the analytical solution are prepared in Sections 4-5. Section 6 is devoted to the applications of the hybrid local meshless method to solve three examples of distributed optimal control problems. The last section is devoted to a brief conclusion.
\section{Necessary optimality conditions}
In this section we present the necessary optimality conditions for stated PDE constrained optimization problems in introduction. We now wish to find the continuous optimality condition for the Lagrangian
\begin{align*}
\mathcal{L}=\frac{1}{2}\int _0^T&\int _{\Omega}(y-y_d)^2d\Omega dt+\frac{\nu}{2	}\int _0^T\int _{\Omega}
&+\int _0^T\int _{\Omega}\left(-\frac{\partial y}{\partial t	}+\Delta y-u\right) p_\Omega d\Omega dt+\int _0^T\int _{\partial\Omega}(y-h)p _{\partial\Omega}dsdt
\end{align*}

where the Lagrange multiplier (adjoint variable) $ p $ has components $ p_\Omega $ and $ p_{\partial\Omega} $ on the interior and boundary of $ \Omega $, respectively. Here the initial condition $ y(x,0)=y_0(x) $ is absoebed into the Lagrangian.\\
From here, the continuous optimality conditions are obtained by differentiating $ \mathcal{L} $ with respect to the adjoint, control, and state variables. Firstly, differentiating with respect to $ p $ returns the forward problem
\begin{align*}
-\dfrac{\partial y}{\partial t}+\triangle y\, & =\,u\,\,\,  in\,\,  Q:=\Omega\times (0,T),\\
y\, & =\,h\,\,\, on\,\,  \Sigma :=\Gamma\times (0,T),\\
y\, & =\,y_0\,\,\,  at\,\,  t=0.
\end{align*} 
Next, differentiating with respect to $ u $ gives us the gradient equation
$$\nu u-p=0  $$
Finally, differentiating with respect to $ y $ gives the adjoint problem
\begin{align*}
\dfrac{\partial p}{\partial t}+\triangle p\, & =\,y_d-y\,\,\,  in\,\,  Q:=\Omega\times (0,T),\\
p\, & =\,0\,\,\, on\,\,  \Sigma :=\Gamma\times (0,T),\\
p\, & =\,0\,\,\,  at\,\,  t=T.
\end{align*}  
We now use the proportionality of control and adjoint, given by the gradient equation, to observe that the conditions reduce to a coupled system of PDEs
\begin{align}
-\dfrac{\partial y}{\partial t}+\triangle y&=\dfrac{1}{\nu}p\\
\dfrac{\partial p}{\partial t}+\Delta p\, & =\,y_d-y
\end{align}

\section{Reproducing kernel space}
In recent years, there is much interest in the use of reproducing kernel for the solution of nonlinear
physical and engineering problems [23-28]. Reproducing kernel theory is used for finding accurate solution of a special class of
nonlinear operator equations by Li and Cui in [22]. Geng and Cui obtained approximate solution for system
of second order nonlinear differential equations by use of exponential kernel in [23]. Jiang and Lin applied
reproducing kernel theory to obtain approximate solution of time-fractional telegraph equation in [24]. Geng
constructed a new RKHS for obtaining convergent series solution of fourth-order two point boundary value
problems in [25]. Arqub et al. used reproducing kernel theory for approximate solution of various type of
Fredholm integro-differential equations in [26]. Bushnaq et al. suggested reproducing kernel method for
fractional fredholm integro-differential equations that obtained approximate solution and its derivatives are
uniformly convergent in [27]. Geng et al. proposed a modified form of reproducing kernel method for solving
singular perturbation problem in [28]. Arqub et al. presented RKHS theory to gain series solution of fuzzy
differential equations in [29]. Recently, periodic boundary value problem of two-point second-order mixed
integro-differential equation is solved in [30]. Some important concepts and some useful reproducing kernel spaces are given in follows.\\
We consider Hilbert spaces over the field of real numbers, $ \mathcal{R} $. Given a set $ X $, we show the set of all functions from $ X $ to $ \mathcal{R} $ with $ \mathbf{F}(X,\mathcal{R}) $. \\
$ \bf{Definition\, 3.1.} $ Given a set X, we will say that $ \mathbf{H} $ is a reproducing kernel Hilbert space (RKHS) on $ X $ over $ \mathcal{R} $, provided that:
\begin{enumerate}
\item $ \mathbf{H} $ is a vector subspace of $ \mathbf{F}(X,\mathcal{R}) $.
\item $ \mathbf{H} $ is endowed with an inner product, $ \langle ,\rangle $, making it into a Hilbert space, 
\item for every $ y\in X $, the linear evaluation functional, $ E_y:\mathbf{H}\rightarrow\mathbf{F} $, defined by $ E_y(f)=f(y) $, is bounded.
\end{enumerate}
If $ \mathbf{H} $ is a RKHS on $ X $, then since every bounded linear functional is given by the inner product with a unique vector in $ \mathbf{H} $, we have that for every $ y\in X $, there exists a unique vector, $ k_y\in \mathbf{H} $, such that for every $ f\in\mathbf{H} $, $ f(y)=$ $\langle f,k_y \rangle $.\\
$ \bf{Definition\, 3.2.} $ The function $ k_y $ is called the reproducing kernel for the point $ y $.\\
The function defined by $$ k(x,y)=k_y(x) $$ is called the reproducing kernel for $\mathbf{H}  $.\\
$ \bf{Definition \,3.3.} $ A Hilbert space $ \mathbf{H} $ of functions on a set $ \Omega $ is called a RKHS if there exists a reproducing kernel $ \mathbf{K} $ of $ \mathbf{H}  $.\\
The existence of the reproducing kernel of a Hilbert space H is due
to the Riesz Representation Theorem. It is known that the reproducing kernel
is unique. We note that it is possible to define several different inner products in the
same class of functions $ \mathbf{H} $, so that $ \mathbf{H} $  is complete with respect to each one
of the corresponding norms. To each one of the Hilbert space $ (H, \langle .,. \rangle) $ there
corresponds one and only one kernel function $ \mathbf{K} $. Therefore, $ \mathbf{K} $ depends not
only on the class of functions in $ \mathbf{H} $, but also on the choice of the inner product
that $ \mathbf{H} $ admits.\\
In order to solve problem (1.1), reproducing kernel spaces $ W_m[a,b] $ with $ m=1,\,2\,3\,... $ are defined in the following, for more details and proofs we refer to [10].\\
 $\bf{Definition\, 3.4.} $ The inner product space $ W_m[a,b] $ is defined as 
\begin{align}
 W_m[a,b]=\lbrace{u(x)\vert u^{'},\,u^{''},\cdot, u^{(m)}\,is\, absolutely\,\, continuous\,\, real\,\, valued \,\,function,\, u^{(m+1)}\in L^2[a,b],\,Bu=0}\rbrace . 
 \end{align}
The inner product in $ W_m[a,b] $ is given by 
\begin{align}
(u(.),v(.))_{W_2}=u^{'}(a)v^{'}(a)+\int _a^b u^{(m+1)}(x)v^{(m+1)}(x)dx,
\end{align}
and the norm $ \Vert u\Vert _{W_m}=\sqrt{(u,u)_{W_m}} $ where $ u,v\in W_m[a,b] $.\\
$ \bf{Theorem\, 3.1.} $ The space $ W_m[a,b] $ is a reproducing kernel space. That is, for any $ u(.)\in W_m[a,b] $ and each fixed $ x\in [a,b] $, there exists $ K(x,.)\in W_m[a,b] $, such that $ (u(.),K(x,.))_{w_m}=u(x) $. The reproducing kernel $ K(x,.) $ can be denoted by 
\begin{align}
k(x,y)=
\begin{cases}
\sum _{i=1}^{2m}c_i(y)x^{i-1}, \,\,\,x\leq y,\\
\sum _{i=1}^{2m}d_i(y)x^{i-1}, \,\,\,x>y.
\end{cases}
\end{align}
According to Definition 3.4, spaces $ W_1, W_2 $ and $ W_1^{'} $ are defined as follows:
\begin{align*}
& W_1[0,T]=\lbrace{u(x)\vert u^{'}\,is\, absolutely\,\, continuous\,\, real\,\, valued \,\,function,\, u^{''}\in L^2[0,T],\,u(0)=0}\rbrace . \\
 &W_1^{'}[0,T]=\lbrace{u(x)\vert u^{'}\,is\, absolutely\,\, continuous\,\, real\,\, valued \,\,function,\, u^{''}\in L^2[0,T],\,u(T)=0}\rbrace . \\
& W_2[a,b]=\lbrace{u(x)\vert u^{'},\,u^{''}\,is\, absolutely\,\, continuous\,\, real\,\, valued \,\,function,\, u^{(3)}\in L^2[a,b],\,u(a)=u(b)=0}\rbrace . 
 \end{align*}
$ \bf{Definition\,3.5.} $ $ W_{(1,2)}(\Omega)=W_1[0,T]\otimes W_2[a,b]=\lbrace u(x,t)\vert\frac{\partial^3 u}{\partial x^2\partial t} $ is completely continuous in $ \Omega,\,\frac{\partial^5u}{\partial x^3\partial t^2}\in L^2(\Omega) $ and $ u(x,0)=u(a,t)=u(b,t)=0\rbrace $. The inner product and the induced norm in $ W(\Omega) $ are defined respectively by 
\begin{align*}
\langle u(x,t),v(x,t)\rangle _W=&\sum _{i=0}^1\int _0^T\left[ \frac{\partial ^2}{\partial t^2}\frac{\partial ^i}{\partial x^i}u(a,t)\frac{\partial ^2}{\partial t^2}\frac{\partial ^i}{\partial x^i}v(a,t)\right]dt+\sum _{j=0}^1\left\langle\frac{\partial ^j}{\partial t^j}u(x,0),\frac{\partial ^j}{\partial t^j}v(x,0)\right\rangle _{W_2}\\
&+\int _0^T\int _a^b\left[ \frac{\partial ^3}{\partial t^3}\frac{\partial ^2}{\partial x^2}u(x,t)\frac{\partial ^3}{\partial t^3}\frac{\partial ^2}{\partial x^2}v(x,t)\right]dt 
\end{align*}
and
\begin{align*}
\Vert u\Vert _W=\sqrt{\langle u,u\rangle _W},\,\, u\in W(\Omega)
\end{align*}
$ \bf{Theorem \,3.2.} $ $ W(\Omega) $ is a reproducing kernel space and its reproducing kernel is 
\begin{align}
k_{(r,s)}^w(x,t)=k_r^{w_1}(x)k_s^{w_2}(t),
\end{align}
such that for any $ u(x,t)\in W(\Omega) $
\begin{align*}
u(r,s)=\langle u(x,t),k_{(r,s)}^w(x,t)\rangle .
\end{align*}
$ \bf{Proof.} $ See [22].\\
The $ W_{(1,2)}^{'} $ space is defined similar to the definition of $ W $ only $ W_1 $ must be replaced with $ W_1^{'} $.
\section{Implement the RKHS based method}
There are three essential steps to complete the solution process of the coupled system of PDEs (2.1)-(2.2) by RKHS 
\begin{enumerate}
\item Construct the reproducing kernel space and corresponding reproducing kernel
satisfying all the conditions for determining solution. The ability of the reproducing
kernel space to absorb all the conditions for determining solution
reflects the feature of the reproducing kernel method.
\item Use the reproducibility of $ \varphi (x)=K_{x_i}(x):\langle f,\varphi _i\rangle _{H_1}=f(x_i) $.
\item Use  $ \psi _i(x)=\Delta \varphi _i(x) $ to construct the base of the space. The solution of the coupled system of PDEs (2.1)-(2.2) is expressed in the form of series.
\end{enumerate}
We first rewrite equation (2.1)-(2.2) and related initial and boundary conditions as follows
\begin{align}
L_1y=F_1(x,t,y,p)+G_1(x,t)\\
L_1p=F_2(x,t,y,p)+G_2(x,t)
\end{align}
with initial conditions
\begin{align}
y(x,0)=p(x,0)=0,\,\,\,x\in [a,b]
\end{align}
and the boundary conditions
\begin{align}
y(a,t)=y(b,t)=0,\,\,\,       t\in[0,T]\\ 
p(a,t)=p(b,t)=0,
\end{align}
where\\
\begin{align*}
&L_1=-\frac{\partial}{\partial t}+\Delta ,\,\,\, L_2=\frac{\partial}{\partial t}+\Delta \\
&F_1(x,t,y,p)=\frac{1}{\nu}p,\,\,\,F_2(x,t,y,p)=y_d-(y+\hat{y})\\
&G_1(x,t)=\frac{\partial Y}{\partial t}-\Delta y_0,\,\,\,G_2(x,t)=0,
\end{align*}
in which
\begin{align*}
Y=\dfrac{x-b}{a-b}h_1(t)+\dfrac{x-a}{b-a}h_2(t),\,\,\,\hat{y}=Y+y_0-Y(x,0).
\end{align*}
Since $ y(x,t) $ and $ p(x,t) $ are sufficiently smooth, $ L_1:W_{(1,2)}\rightarrow \overline{W}_{(1,2)} $ and $ L_2:W^{'}_{(1,2)}\rightarrow \overline{W^{'}}_{(1,2)}  $ are bounded linear operators. \\In the following we discuss the operator $ L_1 $ and Equation 2.1,  for the operator $ L_2 $ and Equation 2.2 is the same.\\
$ \bf{Theorem\,4.1.} $ The operator $ L_1:W_{(1,2)}\rightarrow \overline{W}_{(1,2)} $ is a bounded operator.\\
$ \bf{Proof.} $ Note that 
\begin{align*}
&\Vert (L_1y)(x,t)\Vert ^2=\Vert -u_t+u_{xx}\Vert\leq\Vert u_t\Vert ^2+\Vert u_{xx}\Vert ^2,\\
&y(x,t)=\langle y(r,s),K^{w_1}(r,x)K^{w_2}(s,t)\rangle \\
&y_t=\langle y(r,s),K^{w_1}(r,x)\dfrac{\partial}{\partial t}K^{w_2}(s,t)\rangle \\
&y_{xx}=\langle y(r,s),\dfrac{\partial ^2}{\partial x^2}K^{w_1}(r,x)K^{w_2}(s,t)\rangle \\
&\vert y_t(x,t)\vert\leq \Vert y\Vert \Vert K^{w_1}(r,x)\Vert\Vert \dfrac{\partial}{\partial t}K^{w_2}(s,t)\Vert\\
&\vert y_{xx}(x,t)\vert\leq \Vert y\Vert \Vert \dfrac{\partial ^2}{\partial x^2}K^{w_1}(r,x)\Vert\Vert K^{w_2}(s,t)\Vert\\
\end{align*}
Also note that
$$\Vert K^{w_1}(r,x)\Vert =\sqrt{\langle K^{w_1}(r,x),K^{w_1}(r,x)\rangle}=\sqrt{K^{w_1}(x,x)},\,\,\,\Vert K^{w_2}(s,t)\Vert =\sqrt{K^{w_1}(t,t)} $$
are continuous functions,such that $\Vert K^{w_1}(r,x)\Vert \leq M_1 $, $\Vert K^{w_2}(s,t)\Vert\leq M_2$. Meanwhile, setting $ \Vert \dfrac{\partial ^2}{\partial x^2}\Vert=M_3 $, $\Vert \dfrac{\partial}{\partial t}\Vert=M_4$, we have
$$ \vert y_t(x,t)\vert\leq\Vert y\Vert M_1M_4\vert y_{xx}\vert\leq\Vert y\Vert M_2M_3. $$
Hence $$\Vert (Ly)(x,t)\Vert ^2\leq\Vert y\Vert (M_1^2M_4^2+M_2^2M_3^2)$$
the proof is complete.\\
We choose a countable dense subset $ \lbrace x_j,t_j\rbrace _{j=1}^{\infty} $ in $ Q $, and define 
$$ \varphi _j(x,t)=K^{\overline{w}}_{(x_j,t_j)}(x,t),\,\,\,\,\psi _{j1}(x,t)=L^{*}_1 \varphi _j(x,t). $$ 
where $ L_1^{*} $ is the adjoint operator of $ L_1 $. It can be shown that
$$\psi_{j1}(x,t)=L_{1(y,s)}K^w_{(y,s)}(x,t)\vert _{(y,s)=(x_j,t_j)}.  $$
The subscript $ (y,s) $ by the operator $ L_1 $ indicates that the operator $ L_1 $ applies to the function of $ (y,s). $\\
$ \bf{Theorem \,\,4.2.} $ Suppose that $ \lbrace(x_j,t_j)\rbrace _{j=1}^{\infty} $ is dense in $ Q $, then the analytical solution of (4.1) can be represented as $$y(x,t)=\sum _{j=1}^{\infty} b_{j1}\psi _{j1}(x,t),  $$
where the coefficient $ b_j $ are determined by solving the following semi-infinite system of linear equations $$\bf{A_1b_1=C_1}  $$  in which\\
\begin{align*}
&\bf{A_1}=[L_1\psi _{j1}(x,t)\vert _{(x,t)=(x_i,t_i)}]_{i,j=1,2,\cdots},\,\,\bf{b_1}=[b_{11},b_{21},\cdots]^T,\\
&\bf{C_1}=[F_1(x_1,t_1,y_1,p_1)+G_1(x_1,t_1), F_1(x_2,t_2,y_2,p_2)+G_1(x_2,t_2),\cdots]^T.
\end{align*} 
$ \bf{Proof.} $ Since $ \lbrace  (x_j,t_j)\rbrace _{j=1}^\infty $ is dense in Q, then $ \psi _{j1}(x,t) $ is a complete system in $ W(Q). $ So the analytical solution can be repeesented as (4.6). Since
\begin{align*}
&\langle \psi _{i1}(x,t),\psi _{j1}(x,t)\rangle _W=\langle L^*_1\phi _i(x,t),\psi _{j1}(x,t)\rangle _W=\langle \phi _i(x,t),L_1\psi _{j1}(x,t)\rangle _{\overline{W}}=L_1\psi _{j1}(x,t)\vert _{(x,t)=(x_i,t_i)},\\
&\langle y(x,t),\psi _{i1}(x,t)\rangle _W=\langle y(x,t),L_1^*\phi _i(x,t)\rangle _W=\langle L_1y(x,t),\phi _i(x,t)\rangle _W=F_1(x_i,t_i,y(x_i,t_i),p(x_i,t_i))+G_1(x_i,t_i),
\end{align*} 
According to the best approximation in Hilbert spaces [17], the coefficients $ b_{j1} $ are
determined by solving the semi-infinite system of the linear equations (4.7) and the
proof is complete.\\
The
analytical solution of each equation can be obtained directly from (4.6). In practice, we only need to use a finite sum of Eq. (4.6) to approximate $ y(x,t) $ and $ p(x,t) $. So,the approximate solution of equation is the n-term intercept of analytical solution
which can be determined by solving a $ n\times n $ system of linear equations.
\section{Convergence analysis} 
We assume that $ \lbrace(x_j,t_j)\rbrace _{j=1}^{\infty} $ is dense in $ Q $,. We discuss the convergence of the approximate solutions constructed
in Section 4. \\
$\bf{Theorem\,5.1.}  $ For each $ y(x,t) $, let $ \varepsilon _n^2=\Vert y(x,t)-y_n(x,t)\Vert ^2 $, then sequence $ \lbrace\varepsilon _n\rbrace $ is monotone decreasing and $ \varepsilon\rightarrow 0\,\,(n\rightarrow\infty). $\\
$ \bf{Proof.} $ \\
Since 
\begin{align*}
K_{(y,s)}(x,t)&=\langle K_{(y,s)}(\xi,\eta),K_{(x,t)}(\xi,\eta)\rangle _W\\
&=\langle K_{(x,t)}(\xi,\eta),K_{(y,s)}(\xi,\eta)\rangle _W=K_{(x,t)}(y,s),\\
\Vert K_{(x,t)}(y,s)\Vert _W&=\langle K_{(x,t)}(y,s),K_{(x,t)}(y,s)\rangle _W=K_{(x,t)}(x,t)<\gamma .
\end{align*}
Also sicne $ L_1 $ is bounded thus
$$\Vert\psi _i(x,t)\Vert=M_i .  $$
Now because\\
\begin{align*}
\varepsilon _n^2&=\Vert y(x,t)-y_n(x,t)\Vert ^2 \\
&=\Vert\sum _{i=n+1}^{\infty}\langle y(x,t),\psi _i(x,t)\rangle\psi _i(x,t)\Vert ^2\\
&=\sum _{i=n+1}^{\infty}\langle y(x,t),\psi _i(x,t)\rangle ^2M^2_i ,
\end{align*}
we have
\begin{align*}
\varepsilon _{n-1}^2&=\Vert y(x,t)-y_{n-1}(x,t)\Vert ^2 \\
&=\Vert\sum _{i=n}^{\infty}\langle y(x,t),\psi _i(x,t)\rangle\psi _i(x,t)\Vert ^2\\
&=\sum _{i=n}^{\infty}\langle y(x,t),\psi _i(x,t)\rangle ^2M^2_i .
\end{align*}
Clearly $ \varepsilon _{n-1}\geq\varepsilon _n $. Hence $ \lbrace \varepsilon _n\rbrace $ is monotone decreasing and $ \varepsilon _n\rightarrow 0\,(n\rightarrow\infty) $.
\section{Numerical simulations}
In this section, numerical experiments are conducted to validate the proposed method. We have solved the following three distributed optimal control problems. The following Examples are from \cite{Nazemi}. The simulation
is conducted on Matlab 7. \\

$\bf{Example\, 5.1.}$ In the problem case (2), we take $T=1$, $x\in [0,1]$ and 
\begin{align*}
&y^d=\nu(2(t-1)^3x(x-1)-6t^2(t-1)^2-4t(t-1)^3+12t(t-1)^2x(x-1)+3t^2(2t-20x(x-1)+\\&4t(t-1)^3+6t^2(t-1)^2)
+t^2(1-t)^3x(x-1).
\end{align*}
The exact solution
\begin{align*}
&y=t^2(1-t)^3x(x-1),\\
&p=\nu(2t(t-1)^3x(x-1)-2t^2(t-1)^3+3t^2(t-1)^2x(x-1)).
\end{align*}
The graphs of analytical and estimated solutions of $ y(x,t) $ and $ p(x,t) $ for $ t=0,0.1,0.3,0.5,0.7,0.9,1 $ with $ N=200 $ and $ \nu =10^{-6} $ are respectively, plotted in Fig. 1. In Fig. 2, the error functions $ y-\hat{y} $ and $ p-\hat{p} $ with $ \nu =10^{-6} $ are, respectively, plotted. 
\begin{figure}[ht]
\begin{center}
\includegraphics[width=7cm,height=6cm]{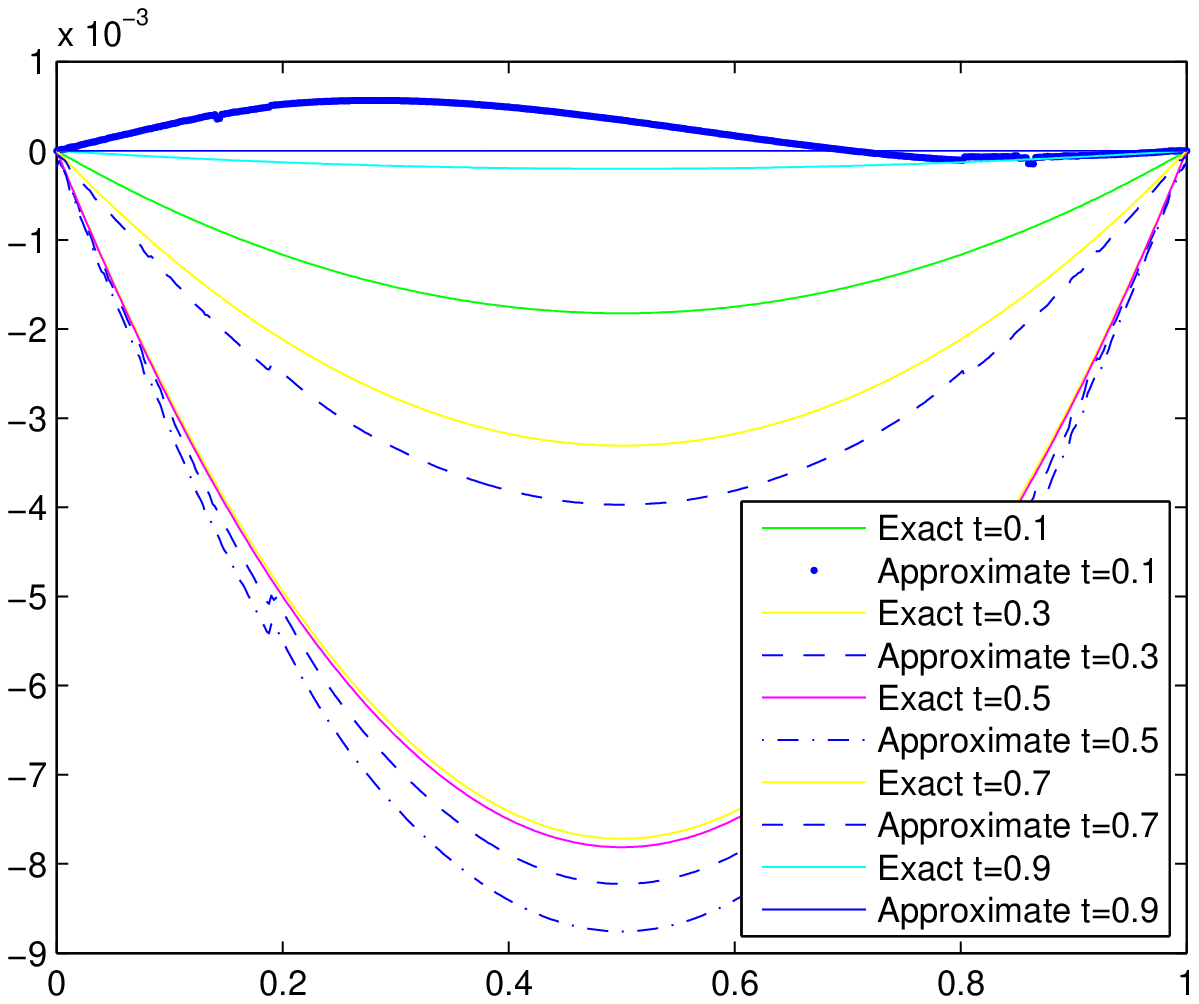}
\includegraphics[width=7cm,height=6cm]{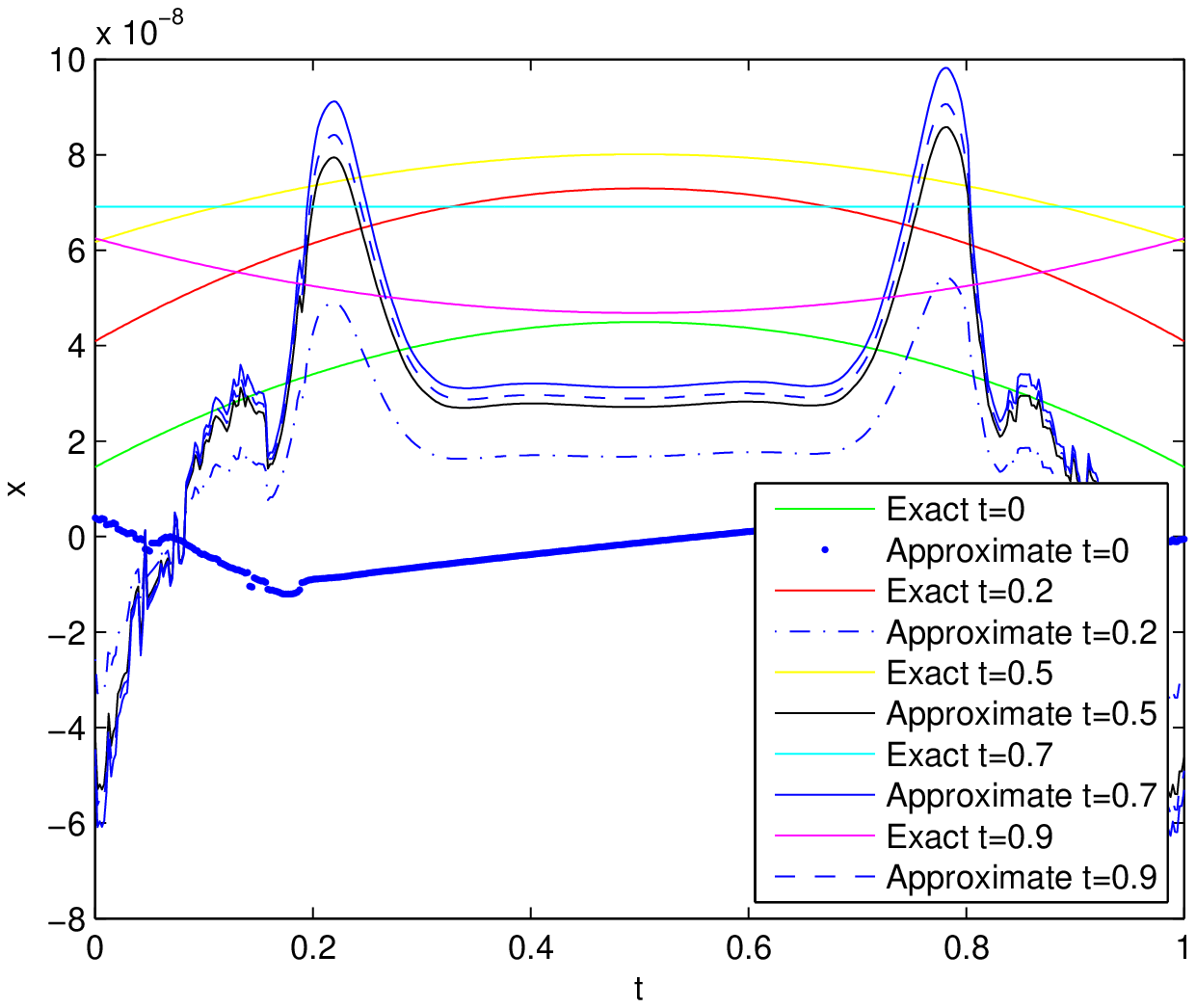}
\vspace{-0.8cm}
\caption{\footnotesize{Comparisons between analytical and approximated solutions of $ y(x,t) $ (left) and $ p(x,t) $ (right) in t=0s, t=0.2s, t=0.5s, t=0.7s, t=0.9s, t=1s with $ \nu =10^{-6} $ in Example 5.1.}}\label{fig1}
\hspace{3mm}
\includegraphics[width=7cm,height=6cm]{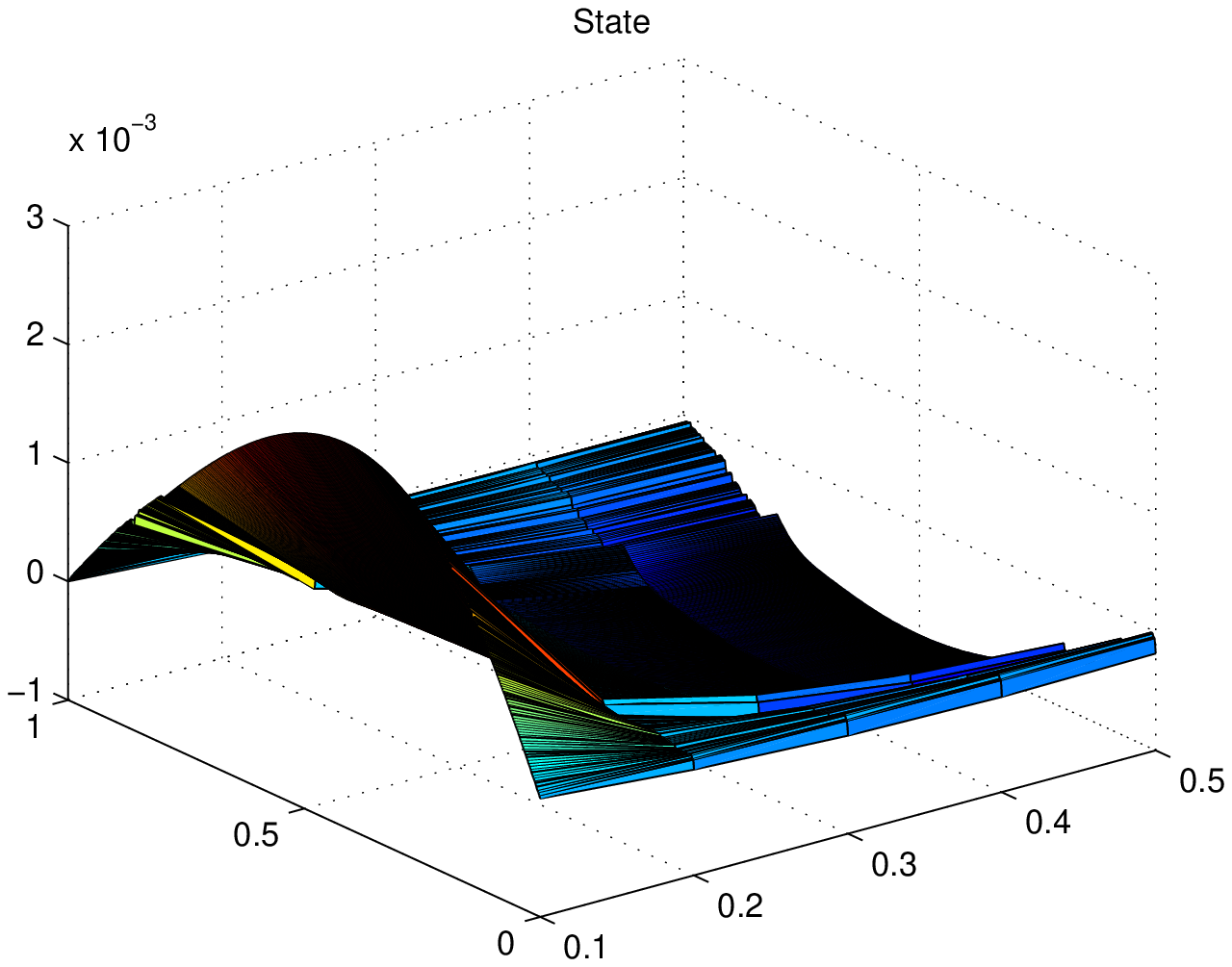}
\includegraphics[width=7cm,height=6cm]{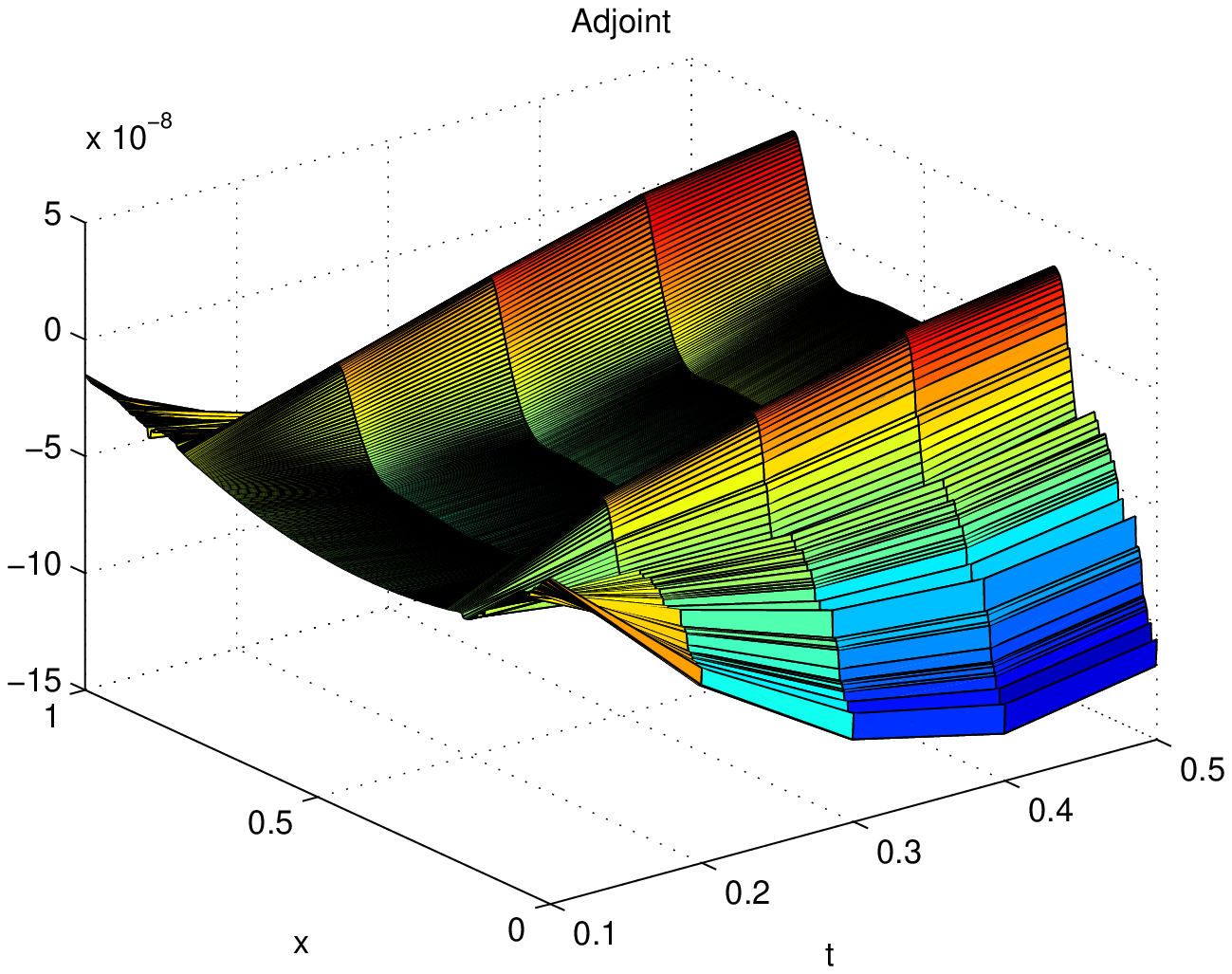}
\vspace{-0.8cm}
\caption{\footnotesize{Plots of $ y-\hat{y} $ and $ p-\hat{p} $  with $ \nu =10^{-6} $ in Example 5.1.}}\label{fig1}
\end{center}
\end{figure}
\\$\bf{Example 5.2.}$ In problem case (2), we take $T=1$, $x\in [0,1]$ and 
\begin{align*}
&y^d=\nu(2\pi ^2t(t-1)^2(t-2)^2+pi ^2t^2(2t-2)(t-2)^2+\pi ^2t^2(2t-4)(t-1)^2+\pi ^4t^2(t-1)^2(t-2)^2-2t^2(t-1)^2\\&-2t^2(t-2)^2-2(t-1)^2(t-2)^2-2t^2(2t-2)(2t-4)-4t(2t-2)(t-2)^2-4t(2t-4)(t-1)^2-\\&2\pi ^2t(t-1)^2(t-2)^2-\pi ^2t^2(2t-2)(t-2)^2-\pi ^2t^2(2t-4)(t-1)^2)sin(\pi x)+t^2(1-t)^2(2-t)^2sin(\pi x).
\end{align*}
The exact solution is 
\begin{align*}
&y=t^2(1-t)^2(2-t)^2sin(\pi x),\\
&p=\nu(-t^2((2t-t)(t-2)^2-t^2(2t-4)(t-1)^2-2t(t-1)^2(t-2)^2-\pi ^2t^2(1-t)^2(2-t)^2sin(\pi x).
\end{align*}
The graphs of analytical and estimated solutions of $ y(x,t) $ and $ p(x,t) $ for $ t=0,0.1,0.3,0.5,0.7,0.9,1 $ with $ N=200 $ and $ \nu =10^{-6} $ are respectively, plotted in Fig. 3. In Fig. 4, the error functions $ y-\hat{y} $ and $ p-\hat{p} $ with $ \nu =10^{-6} $ are, respectively, plotted. 
\begin{figure}[ht]
\begin{center}
\includegraphics[width=7cm,height=6cm]{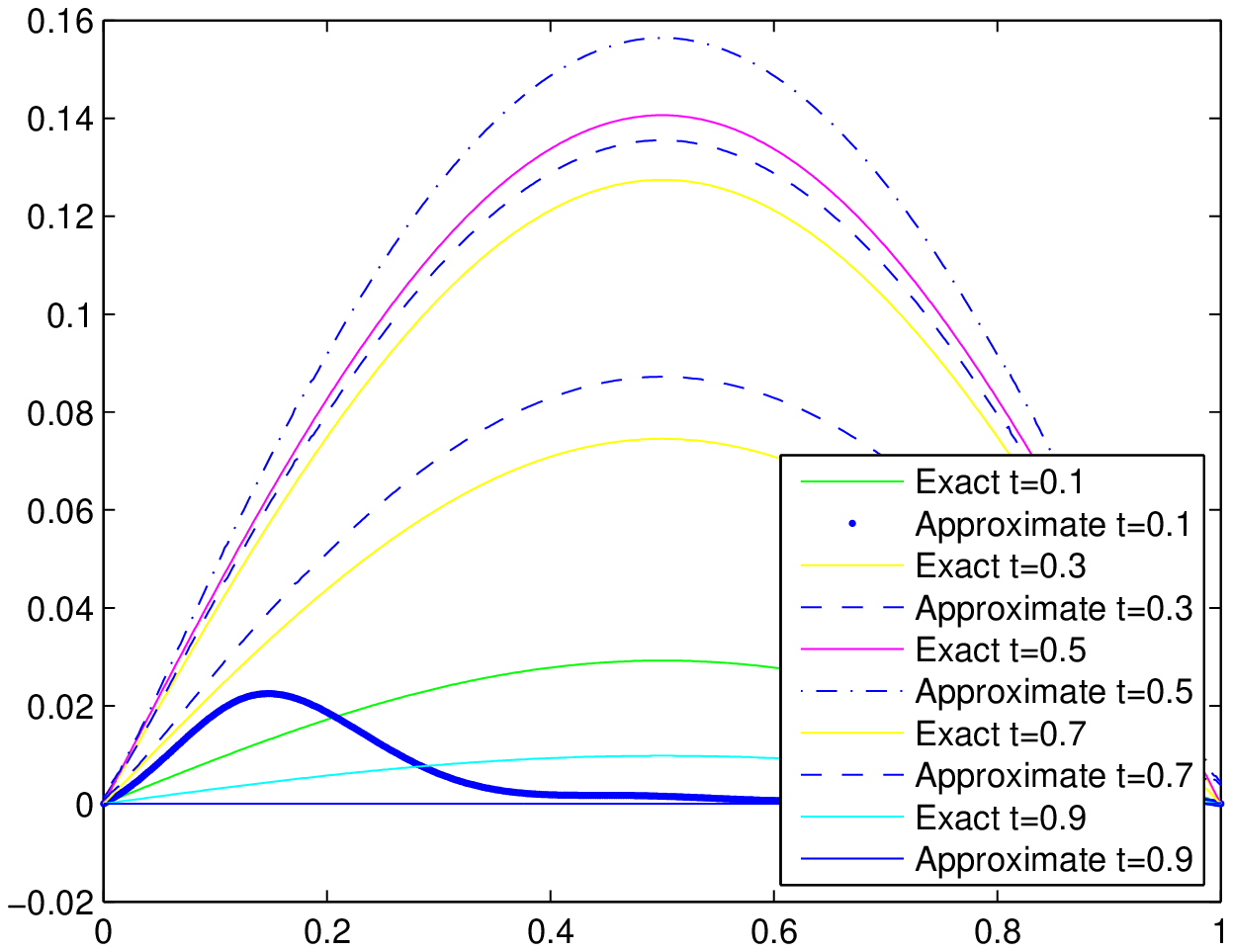}
\includegraphics[width=7cm,height=6cm]{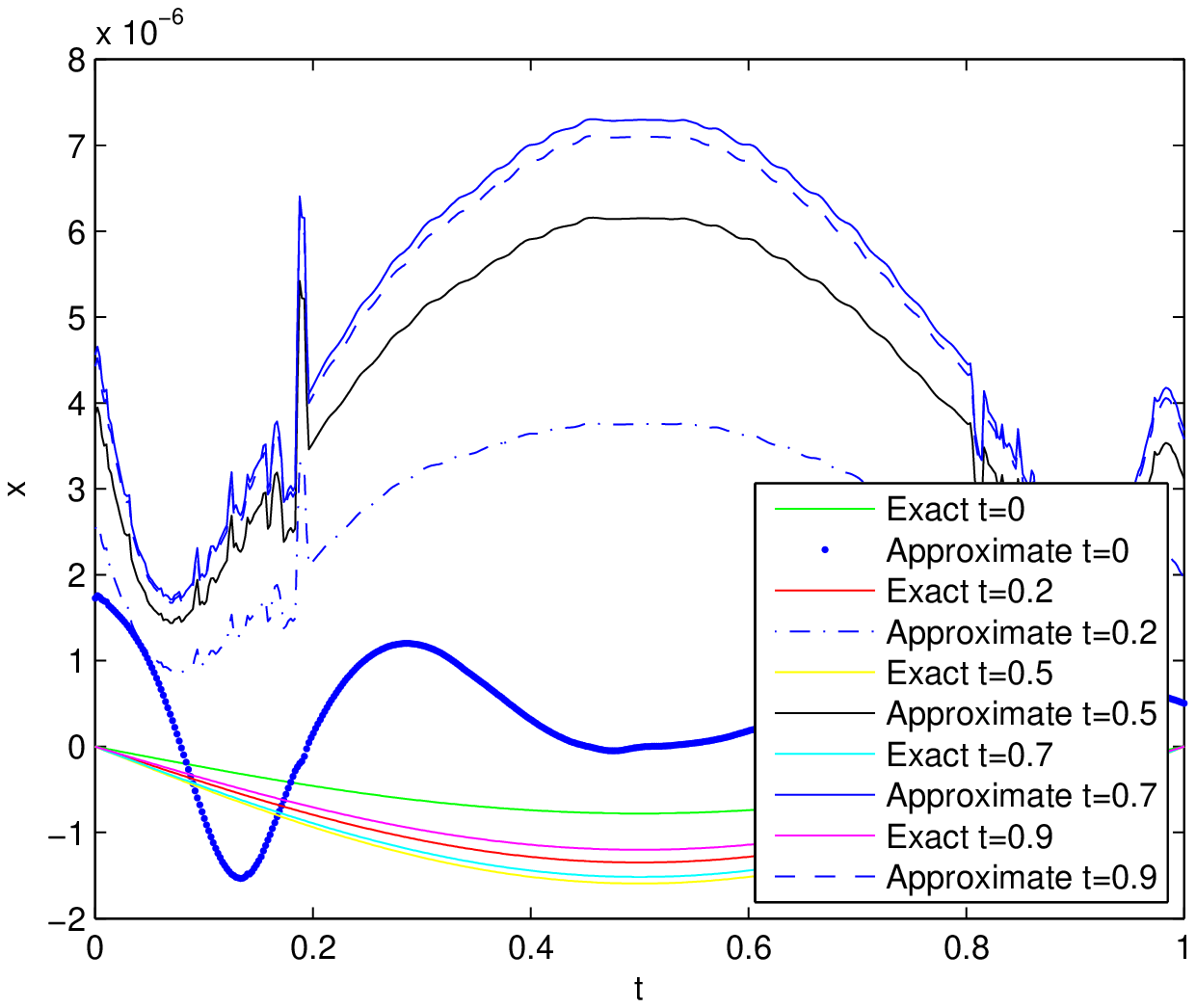}
\vspace{-0.8cm}
\caption{\footnotesize{Comparisons between analytical and approximated solutions of $ y(x,t) $ (left) and $ p(x,t) $ (right) in t=0s, t=0.2s, t=0.5s, t=0.7s, t=0.9s, t=1s with $ \nu =10^{-6} $ in Example 2.}}\label{fig1}
\hspace{3mm}
\includegraphics[width=7cm,height=6cm]{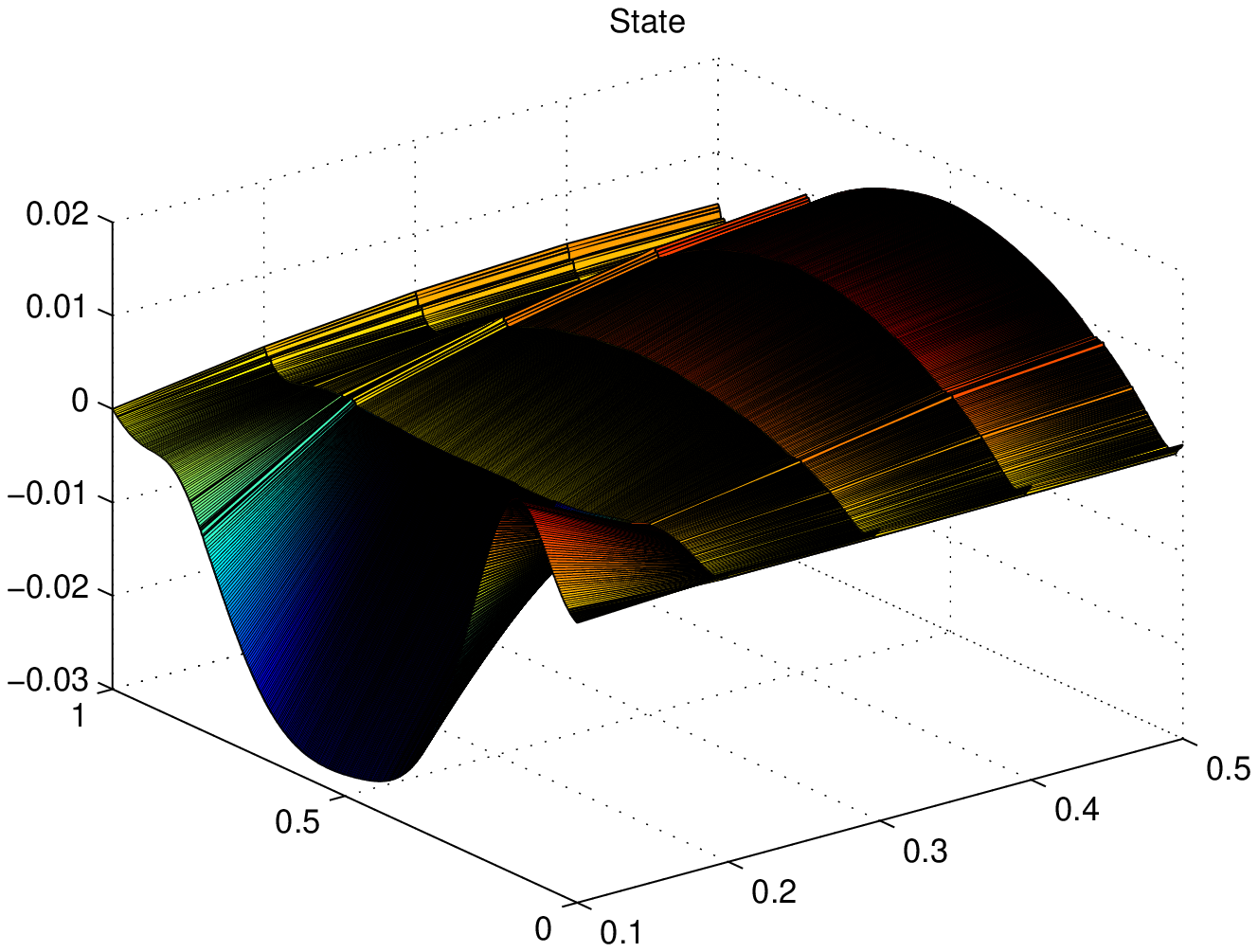}
\includegraphics[width=7cm,height=6cm]{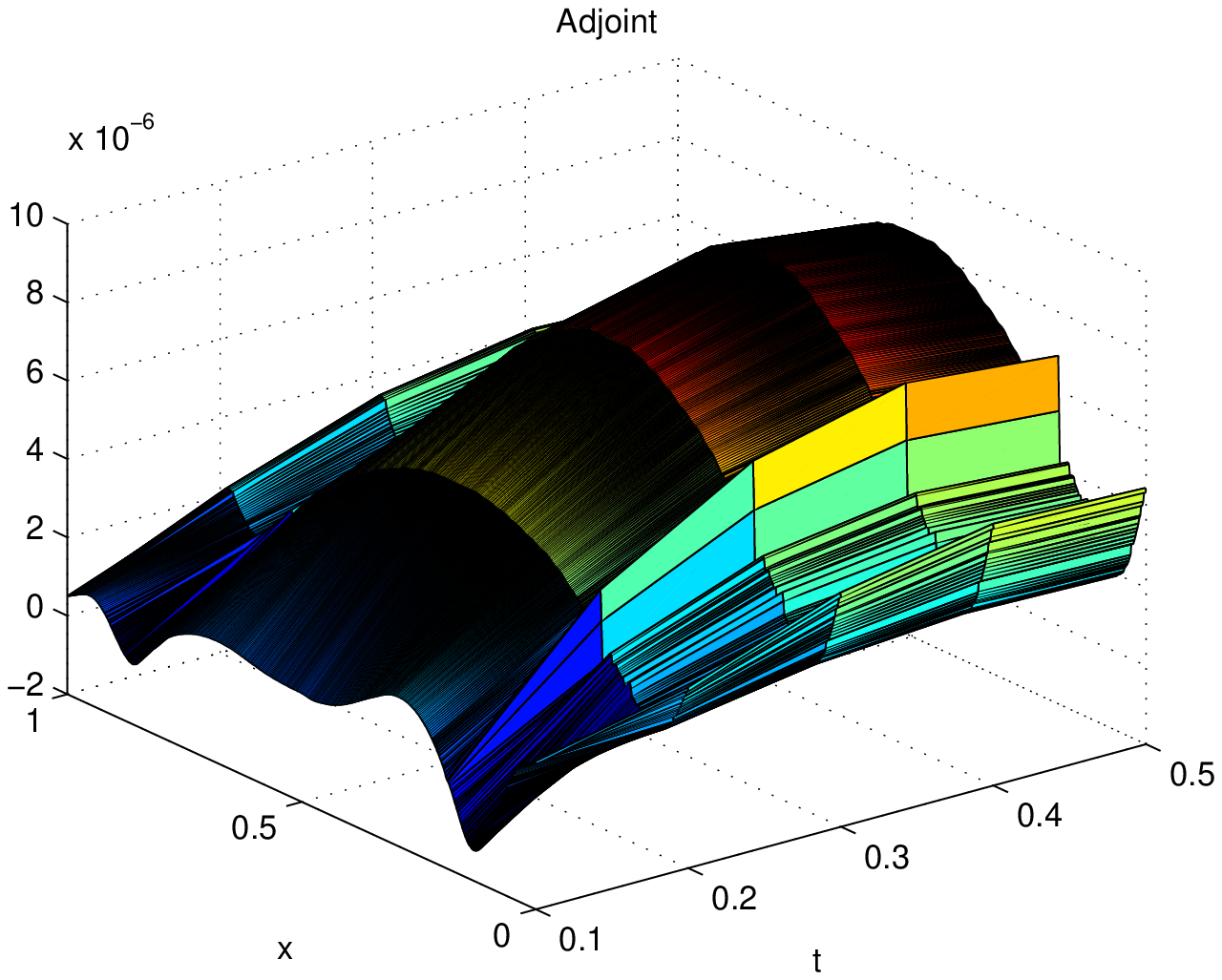}
\vspace{-0.8cm}
\caption{\footnotesize{Plots of $ y-\hat{y} $ and $ p-\hat{p} $  with $ \nu =10^{-6} $ in Example 2.}}\label{fig1}
\end{center}
\end{figure}

$\bf{Example 5.3.}$ In the problem case (2), we take $T=1$, $x\in [0,1]$ and 
\begin{align*}
&y^d=\nu((16\pi ^4t^3(t-1)^3-3t^3(2t-2)-18t^2(t-1)^2-6t(t-1)^3)cos(2\pi x))+\\
&\nu(3t^3(2t-2)+18t^2(t-1)^2+6t(t-1)^3)+t^3(1-t)^3(1-cos(2\pi x)).
\end{align*}
The exact solution
\begin{align*}
&y=t^3(1-t)^3(1-cos(2\pi x)),\\
&p=\nu((-3t^2(t-1)^3-3t^3(t-1)^2)(cos(2\pi x)-1)-4\pi ^2t^3(t-1)^3cos(2\pi x)).
\end{align*}
The graphs of analytical and estimated solutions of $ y(x,t) $ and $ p(x,t) $ for $ t=0,0.1,0.3,0.5,0.7,0.9,1 $ with $ N=200 $ and $ \nu =10^{-6} $ are respectively, plotted in Fig. 5. In Fig. 6, the error functions $ y-\hat{y} $ and $ p-\hat{p} $ with $ \nu =10^{-6} $ are, respectively, plotted. 
\begin{figure}[ht]
\begin{center}
\includegraphics[width=7cm,height=6cm]{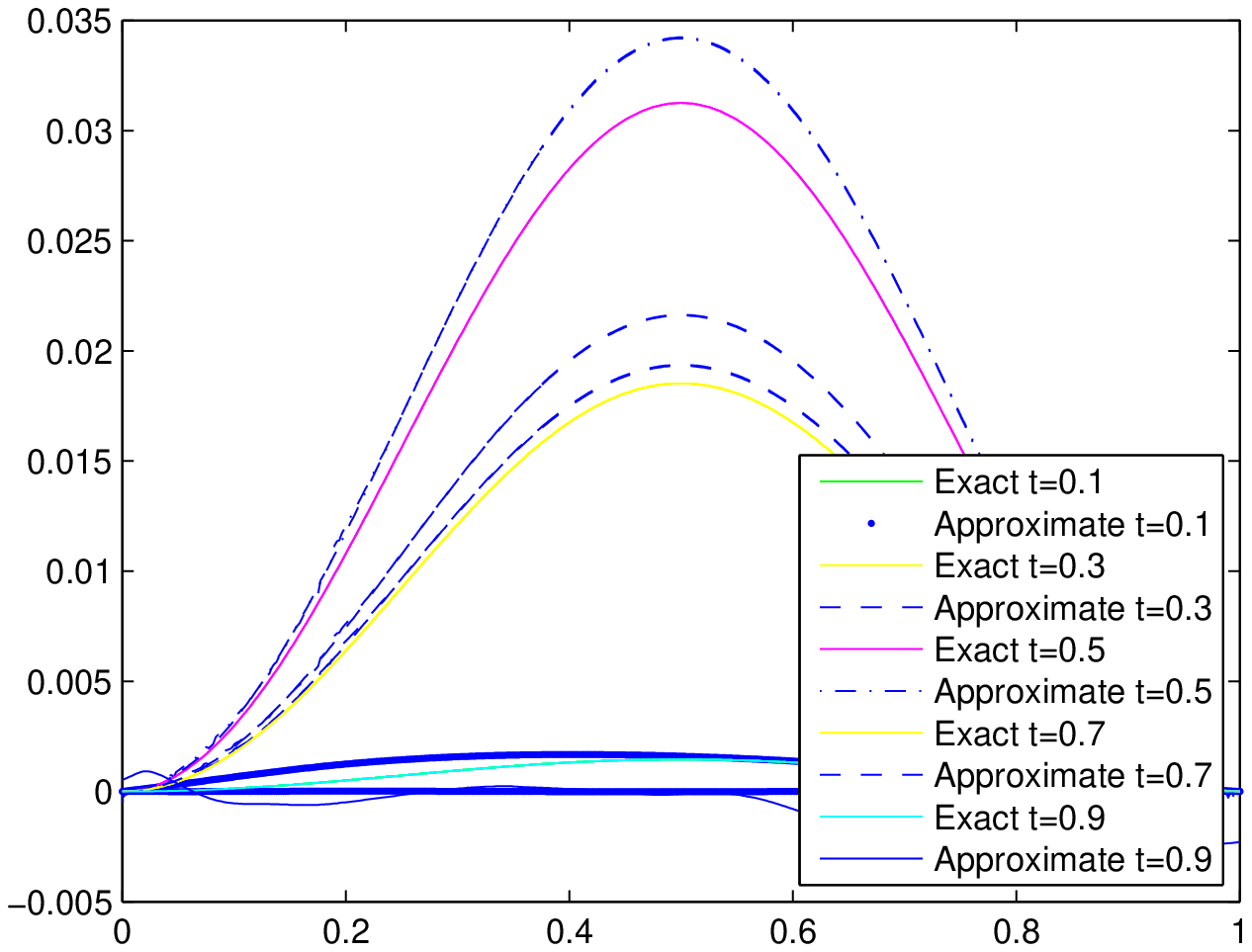}
\includegraphics[width=7cm,height=6cm]{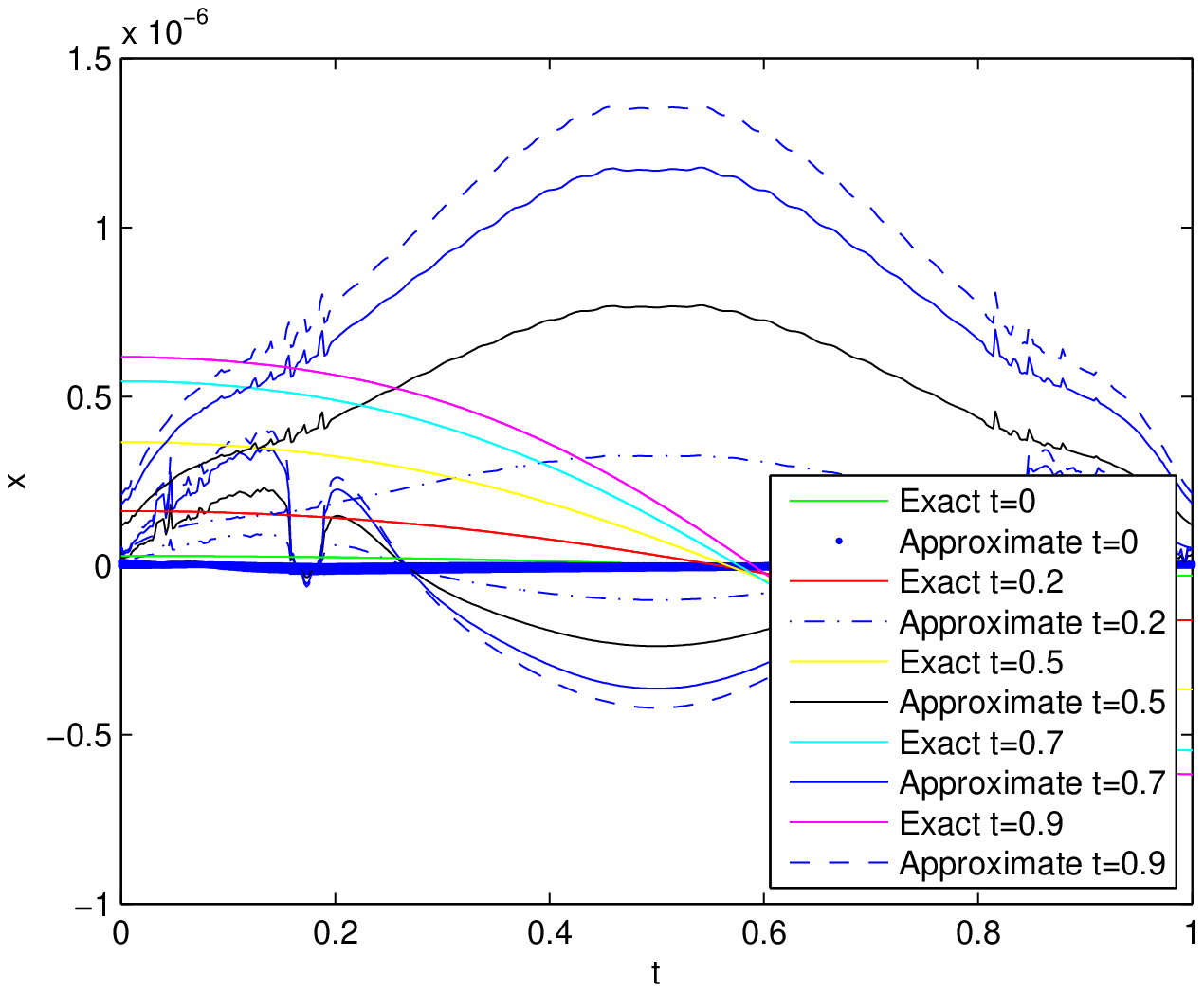}
\vspace{-0.8cm}
\caption{\footnotesize{Comparisons between analytical and approximated solutions of $ y(x,t) $ (left) and $ p(x,t) $ (right) in t=0s, t=0.2s, t=0.5s, t=0.7s, t=0.9s, t=1s with $ \nu =10^{-6} $ in Example 5.3.}}\label{fig1}
\hspace{3mm}
\includegraphics[width=7cm,height=6cm]{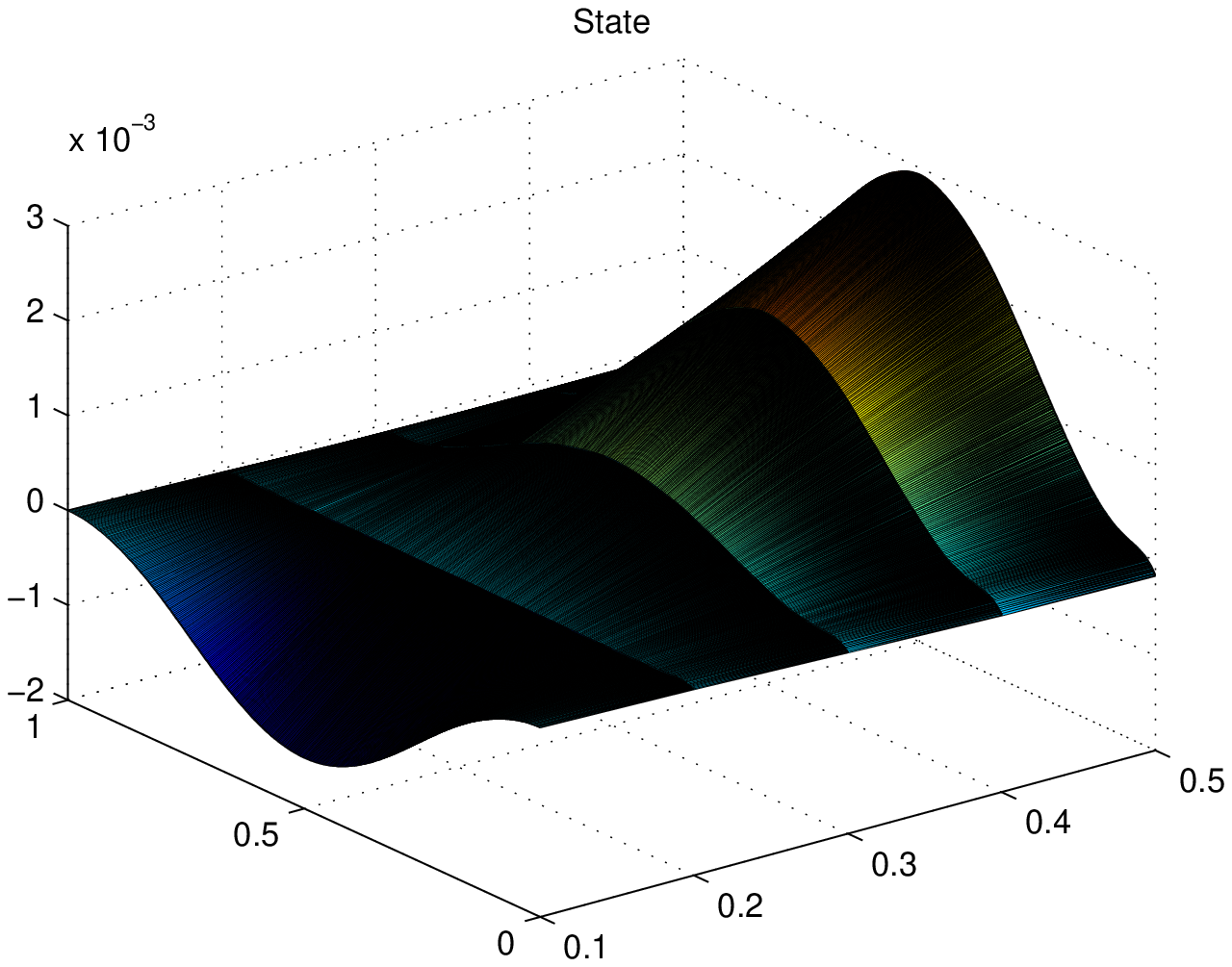}
\includegraphics[width=7cm,height=6cm]{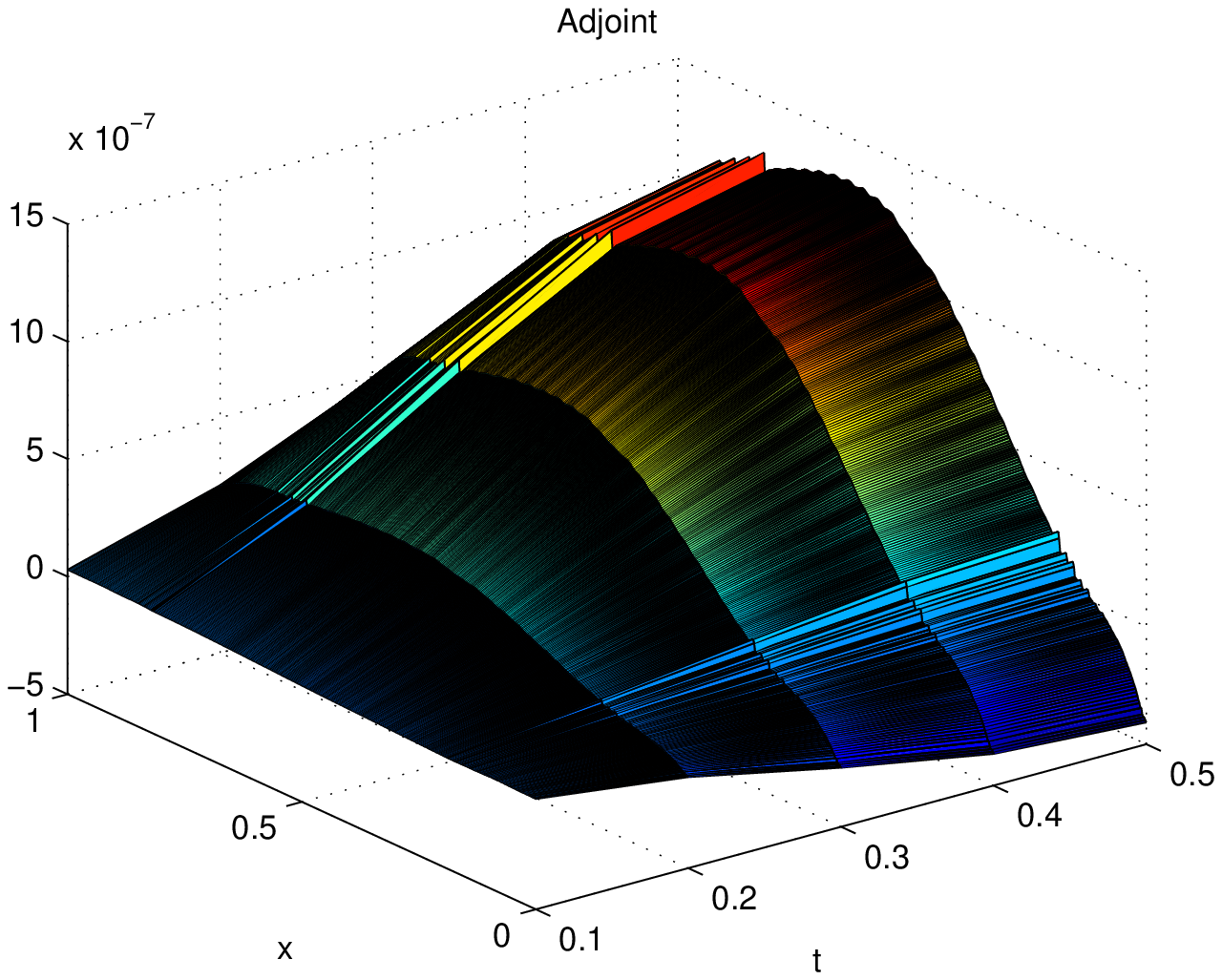}
\vspace{-0.8cm}
\caption{\footnotesize{Plots of $ y-\hat{y} $ and $ p-\hat{p} $  with $ \nu =10^{-6} $ in Example 5.3.}}\label{fig1}
\end{center}
\end{figure}
\section{Conclusion}
In this paper, we use a RPKHS method to solve distributed optimal control problems. The advantages of the used approach lie in the following facts. The method is
mesh free, easily implemented and capable in treating various boundary conditions.
The method needs no time discretization. The used technique is
applied to solve three test problems and the resulting solutions are in good agreement with
the known exact solutions. The numerical results confirmed the efficiency, reliability and accuracy of our method. Furthermore, our method is applicable to more general inverse source problems for parabolic equations, as we will discuss in a forthcoming paper.\\

\end{document}